\documentclass[12pt]{amsart}

\usepackage{amssymb,latexsym,amsmath,amsthm,amscd,verbatim,mathrsfs}

\usepackage[left=3cm,top=2cm,right=3cm,bottom = 2cm]{geometry}

\newtheoremstyle{def}% name
     {10pt}%      Space above
     {10pt}%      Space below
     {}%         Body font
     {}%         Indent amount (empty = no indent, \parindent = para indent)
     {\rmfamily\mdseries\itshape}% Thm head font
     {.}%        Punctuation after thm head
     {.5em}%     Space after thm head: " " = normal interword space;
     {}%         Thm head spec (can be left empty, meaning `normal')
 \theoremstyle{def}
\newtheorem{definition}{Definition}[section]
\newtheorem{remark}[definition]{Remark}

\newtheoremstyle{thm}% name
     {20pt}%      Space above
     {10pt}%      Space below
     {\it}%         Body font
     {}%         Indent amount (empty = no indent, \parindent = para indent)
     {\rmfamily\bfseries\upshape}% Thm head font
     {.}%        Punctuation after thm head
     {.5em}%     Space after thm head: " " = normal interword space;
     {}%         Thm head spec (can be left empty, meaning `normal')

\theoremstyle{thm}
\newtheorem{theorem}[definition]{Theorem}
\newtheorem{lemma}[definition]{Lemma}
\newtheorem{proposition}[definition]{Proposition}
\newtheorem{corollary}[definition]{Corollary}

\DeclareMathOperator{\ord}{ord}

\newcommand{\cA}{\mathcal{A}}

\newcommand{\cE}{\mathcal{E}}
\newcommand{\cF}{\mathcal{F}}

\newcommand{\cH}{\mathcal{H}}
\newcommand{\cI}{\mathcal{I}}

\newcommand{\sM}{\mathscr{M}}

\newcommand{\sX}{\mathscr{X}}

\newcommand{\CC}{\mathbb{C}}

\newcommand{\HH}{\mathbb{H}}

\newcommand{\QQ}{\mathbb{Q}}

\newcommand{\ZZ}{\mathbb{Z}}

\newcommand{\bC}{\mathbb{C}}
\newcommand{\rr}{\epsilon}

\newcommand\blfootnote[1]{%
  \begingroup
  \renewcommand\thefootnote{}\footnote{#1}%
  \addtocounter{footnote}{-1}%
  \endgroup
}

\begin{document}

\title{A geometric perspective on $p$-adic properties of mock modular forms}

\author[L. Candelori]{Luca Candelori}
\author[F. Castella]{Francesc Castella}

\address{Department of Mathematics, LSU, Baton Rouge, LA, 70803, USA}
\email{lcandelori@lsu.edu}

\address{Department of Mathematics, UCLA, Los Angeles, CA, 90095-1555, USA}
\email{castella@math.ucla.edu}

\begin{abstract}
In \cite{BGK}, Bringmann, Guerzhoy and Kane have shown how to correct mock modular forms
by a certain linear combination of the Eichler integral of their shadows
in order to obtain $p$-adic modular forms in the sense of Serre. In this paper,
we give a new proof of their results (for good primes $p$) by employing
the geometric theory of harmonic Maass forms developed by the first author \cite{Candelori}
and the theory of overconvergent modular forms due to Katz and Coleman.
\end{abstract}

\maketitle

\blfootnote{MSC (2010): 11F37, 11F85 }

\section{Introduction}

Over the past decade, there has been a renewed interest in Ramanujan's {\em mock modular forms} and related objects, such as {\em harmonic Maass forms}, whose Fourier coefficients seem to encode interesting arithmetic data not elsewhere found in the classical theory of modular forms. In this article, we offer a new perspective on the $p$-adic properties of the Fourier coefficients of mock modular forms, based on the algebro-geometric theory of $p$-adic modular forms of Katz-Coleman (\cite{katz1}, \cite{coleman2}). Such $p$-adic properties were originally discovered in \cite{BGK}, \cite{GuerzhoyKentOno}, but we believe our perspective simplifies some of the arguments and provides a theoretical platform for further exploration.

In order to state our results precisely, let $\tau = u +  iv \in \mathfrak{h}$, let $\Gamma_0(N)$ be the congruence subgroup of $\mathrm{SL}_2(\mathbb{Z})$ of matrices that become upper-triangular modulo $N$, and let $\chi$ be a Dirichlet character modulo $N$. Denote by $\mathcal{H}_k(\Gamma_0(N),\chi)$ the vector space of all weight $k$ harmonic Maass forms on $\Gamma_0(N)$ and character $\chi$ (see e.g. \cite[\S{2}]{BGK} for definitions). Any harmonic Maass form $F$ has a decomposition
$$
F = F^{+} + F^{-}
$$
into a holomorphic part $F^{+}$ (with poles supported at the cusps) and an anti-holomorphic part $F^{-}$. The function $F^{+}:\mathfrak{h}\rightarrow \CC$ is what is called a {\em mock} modular form, since it does not transform like a modular form, but its Fourier coefficients resemble those of a modular form. Harmonic Maass forms map into spaces of classical modular forms via differential operators. In particular, let $M^{!}_k(\Gamma_0(N),\chi)$ (resp. $S_k(\Gamma_0(N),\chi)$) be the space of weakly holomorphic modular forms (resp. cusp forms) of weight $k$, level $\Gamma_0(N)$ and character $\chi$. If we let
\begin{equation}
\xi_k := 2iv^k\overline{\frac{\partial}{\partial\overline{\tau}}},
\label{eq:xiOperator}
\end{equation}
then $\xi_{2-k}(F) = f \in S_k(\Gamma_0(N),\chi)$ for all $F\in \mathcal{H}_{2-k}(\Gamma_0(N),\overline{\chi})$,
and the resulting cusp form is called the {\em shadow} of $F$. A fundamental question in the subject is to relate the coefficients of 
a mock modular form $F^+$ to the coefficients of its shadow $f$. In order to obtain results in this direction, we first have to restrict to normalized newforms $f$ and then slightly refine the definition of a harmonic Maass form, as follows. Let $K\subseteq \CC$ be a subfield. For $\Gamma\subseteq \mathrm{SL}_2(\ZZ)$ a congruence subgroup, we denote by $S_k(\Gamma, K)$ the space of cusp forms of weight $k$ and level $\Gamma$ whose $q$-expansion coefficients all lie in the field $K$. Let also $M^{!}_k(\Gamma, K)$ be the space of weakly holomorphic modular forms with coefficients in $K$.

\begin{definition}
Let $f\in S_k(\Gamma_1(N),K)$ be a newform defined over $K$.   
A harmonic Maass form $F\in\mathcal{H}_{2-k}(\Gamma_1(N))$ is {\em good} for $f$ if
\begin{itemize}
\item[(i)] The principal parts of $F$ all lie in $K$.
\item[(ii)] $\xi_{2-k}(F) = f/\|f\|^2$, where $\|f\|$ is the Petersson norm of $f$.
\end{itemize}
\label{def:goodHWMF}
\end{definition}

Suppose that $f = \sum_{n=1}^\infty a_nq^n$ is a (normalized) newform as above,  
%(e.g. $K = K_f$, the field generated by the coefficients of $f$). 
let $F$ be a harmonic Maass form that is good for $f$, and write $F = F^{+} + F^{-}$ 
for its holomorphic and anti-holomorphic parts, with
$$
F^{+} = \sum_{n\gg -\infty} c^{+}(n)q^n.
$$

Let $E_f = \sum_{n=1}^\infty n^{1-k}a_nq^n$ be the Eichler integral of $f$, so that 
$D^{k-1}(E_f)=f$, where $D^{k-1}$ is the differential operator on modular forms acting as $(q d/dq)^{k-1}$
on $q$-expansions. It is shown in \cite{GuerzhoyKentOno} (and also in Theorem~\ref{thm:GKO-1} of this paper, by different methods) 
that for any $\alpha \in \CC$ such that $\alpha - c^{+}(1) \in K$, the coefficients of
$$
\cF_{\alpha}:= F^{+} - \alpha E_f
$$
all lie in $K$, so it makes sense to study their $p$-adic properties. 
To this end, let $p\nmid N$ be a prime, fix once and for all a choice of complex and $p$-adic embeddings 
$\overline{\mathbb{Q}}\hookrightarrow\mathbb{C}$ and $\overline{\mathbb{Q}}\hookrightarrow\mathbb{C}_p$, 
and fix a valuation $v_p$ on $\CC_p$ extending the $p$-adic valuation of $\mathbb{Q}$. 
Suppose the Hecke polynomial $T^2 - a_pT + \chi(p)p^{k-1}$ has roots $\beta$, $\beta'$ with $v_p(\beta)\leqslant v_p(\beta')$, and let 
$V$ be the operator acting as $q\mapsto q^p$ on $q$-expansions. The {\em $p$-stabilizations} of $f$ are the $p$-adic modular forms
$$
f_{\beta} := f - \beta' V(f), \quad
f_{\beta'} := f - \beta V(f),
$$
which are easily seen to be eigenvectors for $U$ with eigenvalues $\beta$ and $\beta'$, respectively. 
Here, $U$ is defined by $U(\sum_n a_n q^n) = \sum_n a_{pn} q^n$, and our first main
result shows that, for most values of $\alpha$, the $p$-stabilized shadow $f_{\beta'}$ can be recovered $p$-adically 
from the corrected mock modular form $\mathcal{F}_\alpha$ by an iterated application of the $U$-operator.  

\begin{theorem}[\cite{GuerzhoyKentOno}, Theorem~1.2(i)]\label{introthm}
Assume that $v_p(\beta)\neq v_p(\beta')$ and assume that $v_p(\beta')\neq k-1$.
Then for all but at most one choice of $\alpha$ with $\alpha-c^+(1)\in K$, we have
\[
\lim_{w\to+\infty}\frac{U^wD^{k-1}(\mathcal{F}_\alpha)}{c_\alpha(p^w)}=f_{\beta}.
\]
\end{theorem}

In Section \ref{section:recovering the shadow}, we give a new proof of this result by viewing $f_{\beta}$ and $f_{\beta'}$ as {\em overconvergent modular forms}, in the sense of \cite{coleman2}. Based on ideas of \cite{bdp}, we prove (Theorem \ref{theorem:canonicalBasis} below) that these two modular forms are $p$-adic representatives of cohomology classes in the $f$-isotypical component of a certain parabolic cohomology group attached to the modular curve $X_1(N)$.
%The problem of when these two classes actually span the whole 2-dimensional isotypical component can be completely reduced to the variational Hodge %conjecture, by work of Emerton (\cite{emerton-note}). 
Under the assumptions of Theorem~\ref{introthm}, the classes $f_{\beta}, f_{\beta'}$ form a basis for this space, and so 
the modular form $D^{k-1}(\mathcal{F}_\alpha)$ (which gives a class in the same space, as shown in \cite{Candelori}) 
can be expressed as a linear combination $f_{\beta}$ and $f_{\beta'}$. Our proof of Theorem~\ref{introthm} then follows by analyzing the action of $U$ in cohomology.

This new proof-template can be applied to similar questions in the theory of mock modular forms. For example, in Section \ref{section:mock modular forms as } 
we interpret the exceptional value of $\alpha$ in Theorem~\ref{introthm} as giving the precise value for which $\cF_{\alpha}$ can be $p$-adically `completed' to obtain a $p$-adic modular form.
%, much in the same way $F^{+}$ can be complex-analytically completed to a smooth modular form (i.e. a harmonic Maass form). 
This was initially discovered by Bringmann, Guerzhoy and Kane, and we reprove here their results \cite{BGK} using our $p$-adic analytic/geometric methods. 
Finally, in Section \ref{section:CM} we discuss the case of when $f$ has CM (also considered \cite{BGK} and \cite{GuerzhoyKentOno}), which requires a different treatment
due to the failure of the assumptions in Theorem~\ref{introthm}.

%\noindent\emph{Acknowledgements.} 
We would like to sincerely thank our doctoral adviser Henri Darmon, who gave this project to one of us during our graduate studies at McGill: the cohomological approach to mock modular forms is essentially due to him. We would also like to thank Matt Boylan and Pavel Guerzhoy for their comments on an earlier version of this paper.

\section{Harmonic Maass forms: the geometric point of view}
\label{section:harmonicMaassForms}

We begin by quickly recalling the geometric interpretation of harmonic Maass forms given in \cite{Candelori}, which will be needed in later sections. For $N>4$, the moduli functor $\sM_1(N)$ of generalized elliptic curves with a point of order $N$ is represented by a smooth and proper scheme over $\ZZ[1/N]$. Let $\cE^{\rm{gen}}\rightarrow \sM_1(N)$ be the universal generalized elliptic curve, and let $\underline{\omega}$ be its relative dualizing (invertible) sheaf. Let $X:=  \sM_1(N)\times_{\ZZ[1/N]} \QQ$ and $Y:= X\smallsetminus C$, where $C$ is the cuspidal subscheme, whose ideal sheaf we denote by $\cI_C$. For any subfield $K\subseteq \CC$, we denote by $X_K, Y_K$ the base-change to $K$. We have well-known canonical isomorphisms
$$
M^{!}_k(\Gamma_1(N), K)\simeq H^0(Y_K, \underline{\omega}^k),  \quad
S_k(\Gamma_1(N), K) \simeq H^0(X_K, \underline{\omega}^k\otimes \cI_C), $$
where a modular form $f$ of weight $k$ is identified with the differential $f(dq/q)^k$. Let $\pi:\cE\rightarrow Y$ be the universal elliptic curve with $\Gamma_1(N)$-level structure. The relative de Rham cohomology of $\pi:\cE\rightarrow Y$ canonically extends to a rank 2 vector bundle $\cH^1_{\rm{dR}}$ over $X$. For any $r\geqslant 0$ let
$$
\cH_r:= \mathrm{Sym}^r(\cH^1_{\mathrm{dR}}),
$$
which is a vector bundle of rank $r+1$ over $X$. The Gauss-Manin connection of $\pi:\cE\rightarrow Y$ extends to a connection with logarithmic poles
$
\nabla: \cH^1_{\mathrm{dR}} \rightarrow \cH^1_{\mathrm{dR}}\otimes \Omega^1_X(\log C)
$
over $X$, and the $r$-th symmetric power of $\nabla$ is a connection with logarithmic poles
$$
\nabla_r: \cH_r \longrightarrow \cH_r\otimes \Omega^1_X(\log C).
$$
Define
$$
\mathbb{H}^1_{\rm{par}}(X,\nabla_r):= \mathbb{H}^1(\cH_r\otimes\cI_C \stackrel{\nabla_r}\longrightarrow \cH_r\otimes \Omega^1_X),
$$
where $\mathbb{H}^{\bullet}$ denotes hypercohomology. Over $\CC$, this group is canonically isomorphic to the classical weight $r$ parabolic cohomology obtained by taking periods of cusp forms. The formation of this cohomology group is compatible under base-change by a field extension $K\supseteq \QQ$ and for all such $K$ and $k\geqslant 2$ there is a filtration (\cite[Thm.~2.7.(i)]{Scholl})
$$
0 \longrightarrow H^0(X_K,\underline{\omega}^{k}\otimes\cI_C) \longrightarrow\mathbb{H}^1_{\rm{par}}(X_K,\nabla_{k-2}) \longrightarrow H^1(X_K,\underline{\omega}^{2-k})\longrightarrow 0
$$
of $K$-vector spaces, so that $S_{k}(\Gamma_1(N),K)$ is naturally a subspace of parabolic cohomology. More generally, all parabolic cohomology classes can be represented in terms of classical modular forms. To state this result, recall that for $k\geqslant 2$ there is an %well-known 
algebraic differential operator of order $k-1$:
$$
D^{k-1}: M_{2-k}^{!}(\Gamma_1(N),K) \longrightarrow M_{k}^{!}(\Gamma_1(N),K)
$$
which acts as $\left(q\,d/dq\right)^{k-1} = \left(\frac{1}{2\pi i} \frac{\partial}{\partial\tau}\right)^{k-1} $.

\begin{theorem}[{\cite[Thm.~6]{Candelori}}]
\label{thm:EichlerShimuraStyle}
Let $K\subseteq \CC$ be a subfield and let $S_{k}^{!}(\Gamma_1(N),K)$ %$k\geqslant 2$, 
be the subspace of those modular forms in $M_k^{!}(\Gamma_1(N), K)$ with vanishing constant coefficient in their $q$-expansions at the cusps. 
Then there is a canonical isomorphism:
$$
\mathbb{H}_{\rm{par}}^1(X_K,\nabla_{k-2}) \simeq \frac{S_{k}^{!}(\Gamma_1(N),K)}{D^{k-1}M_{2-k}^{!}(\Gamma_1(N),K)}.
$$
\end{theorem}

Let now $f \in S_k(\Gamma_1(N),K)$ be a newform. Let $\mathbb{H}_{\rm{par}}^1(X_K,\nabla_{k-2})_f$ be the $f$-isotypical component,
%which is a $K$-vector space of dimension 2,
and let
$$
[\phi] \in \mathbb{H}_{\rm{par}}^1(X_K,\nabla_{k-2})_f
$$
be a class represented by an element $\phi \in S_{k}^{!}(\Gamma_1(N),K)$. By the %well-known 
Shimura isomorphism  $\mathbb{H}_{\rm{par}}^1(X_{\CC},\nabla_{k-2})_f \simeq \CC\, f \oplus \CC\, \bar{f}$, we may write
\begin{equation}\label{eq:complex}
[\phi] = s_1[f] + s_2[\bar{f}]% \in \mathbb{H}_{\rm{par}}^1(X_{\CC},\nabla_{k-2})_f,
\end{equation}
for some $s_1,s_2 \in \CC$. Let now $C^{\infty}_Y$ (resp. $\cA^1_Y$) be the sheaf of smooth functions (resp. smooth differential forms) on $Y_{\CC}$.
The differential $\phi - s_1f - s_2\bar{f}$ is smooth over $Y_{\CC}$, and it defines a class in
$$
\HH^1(\cH_{k-2}\otimes C_Y^{\infty}\xrightarrow{\nabla_{k-2}}\cH_{k-2}\otimes \cA_Y^{1}) = \frac{H^0(Y_{\CC},\cH_{k-2}\otimes\cA^1_Y)}{ \nabla_{k-2} H^0(Y_{\CC},\cH_{k-2}\otimes C_Y^{\infty})}.
$$
This class is trivial by construction, and so there exists a smooth $\cH_{k-2}$-valued modular form $\mathbf{F}$  such that
$
\nabla_{k-2}(\mathbf{F}) = \phi - s_1f - s_2\bar{f}.
$
The vector bundle $\cH_{k-2}$ decomposes into line bundles as $\cH_{k-2} \simeq \underline{\omega}^{2-k}\oplus\underline{\omega}^{4-k}\oplus\ldots\oplus \underline{\omega}^{k-2}$, and we let $F := F_{2-k}$ be the component of $\mathbf{F}$ of weight $2-k$. As shown in \cite[Prop.~4]{Candelori}, $F$ is a harmonic Maass form. If we write $F = F^{+} + F^{-}$ for the holomorphic and anti-holomorphic parts of $F$, then
$$
D^{k-1}(F^{+}) = \phi - s_1\,f, \quad \frac{2i\, v^{2-k}}{(-4\pi )^{k-1}}\frac{\partial}{\partial\overline{\tau}}(F^{-}) = s_2\, \bar{f}.
$$
To obtain a `true' harmonic Maass form we should insist that $\phi \notin S_k(\Gamma_1(N),K)$, i.e. $s_2\neq 0$ in Equation \eqref{eq:complex}.
Then we may rescale $\phi$ so that $\langle \phi, f \rangle = 1$ (cup-product), which amounts to letting $s_2 = 1/\langle \bar{f}, f \rangle = 1/(-4\pi)^{k-1}\|f\|^2$. With this choice, it is clear from the above that $\xi_{2-k}(F) = f/\|f\|^2$, so that $F$ is good for $f$ in the sense of Definition~\ref{def:goodHWMF}.

\section{Overconvergent modular forms}

Let $p\geqslant 5$ be a prime and let $\mathbb{C}_p$ be the completion of the algebraic closure of $\mathbb{Q}_p$. We fix a valuation $v_p$ on $\mathbb{C}_p$ such that $v_p(p)=1$ and an absolute value $|\cdot|$ on $\mathbb{C}_p$ which is compatible with $v_p$. Let $K_p\subseteq \CC_p$ be a complete discretely-valued subfield and let $R_p$ be its ring of integers. Suppose $(p,N)=1$, and let $\sX:= \sM_1(N)\times_{\ZZ[1/N]} R_p$ be the base-change to $R_p$. Let
$
E_{p-1} \in H^0(\sX\times_{R_p} K_p,\underline{\omega}^{p-1})
$
be the global section given by the Eisenstein series of weight $p-1$ and level 1, normalized so that its constant coefficient is 1. As shown in \cite[\S{1}]{coleman2}, for any $\rr\in |R_p|$ there is a unique rigid analytic space $X_\rr$ with the property that
$$
X_\rr^{\rm cl} = \{ x\in (\sX\times_{R_p} K_p)^{\rm{cl}}: |E_{p-1}(x)| \geqslant \rr \},
$$
where by the superscript `cl' we have denoted the set of closed points, and also a unique rigid analytic space $X_{(\rr)}$ with the property that
$$
X_{(\rr)}^{\rm cl} = \{ x\in (\sX\times_{R_p} K_p)^{\rm{cl}}: |E_{p-1}(x)| > \rr \}.
$$

When $\rr=1$, the rigid analytic space $X^{\rm{ord}}:= X_1$ is called the {\em ordinary locus} of $X$, since every geometric point of $X^{\rm{ord}}$ reduces mod $p$ to a point classifying an ordinary elliptic curve. The rigid analytic spaces $X_{(\rr)}$, for $0<\rr<1$, can be viewed as `complements of closed disks' and are called {\em open neighborhoods} of $X^{\mathrm{ord}}$. They are examples of {\em wide open spaces}. For all $\rr\in |R_p|$, we have inclusions $X_1 \subseteq X_\rr \subseteq X_{(\rr)} \subseteq X$. The invertible sheaves $\underline{\omega}^k$, for $k\in \ZZ$, restrict to rigid analytic line bundles over $X_\rr$.

\begin{definition}
Let $\rr\in |R_p|$. An {\em overconvergent modular form} of weight $k\in\mathbb{Z}$ is a rigid analytic section
$f \in H^0(X_{(\rr)}, \underline{\omega}^k)$, for $\rr<1$.
\end{definition}

Note that for $\rr=1$ the sections of $\underline{\omega}^k$ over $X^{\rm{ord}}$ are Serre's $p$-adic modular forms of integral weight $k$. Overconvergent modular forms can thus be viewed as  $p$-adic modular forms which converge not just over $X^{\rm{ord}}$ but on a slightly larger  neighborhood of it.

%\begin{comment}
Since $|E_{p-1}(c)| = 1$ at all cusps $c\in C$, we have that $C\subseteq X_{(\rr)}$ for all $\rr\in |R_p|$. Let
$$
Y^{\ord}:= X^{\rm{ord}}\smallsetminus C, \quad Y_\rr:= X_\rr \smallsetminus C, \quad Y_{(\rr)}:= X_{(\rr)} \smallsetminus C
$$
be the rigid analytic spaces obtained by removing the cusps.

%\begin{definition}
%Let $\rr\in |R_p|$. A {\em weakly overconvergent modular form} of weight $k\in\mathbb{Z}$ is a rigid analytic section
%$f \in H^0(Y_{(\rr)}, \underline{\omega}^k)$, for some $\rr<1$.
%\end{definition}

\begin{remark}
For $\rr=1$, sections of  $H^0(Y_{\rr},\underline{\omega}^k)$ correspond to the $p$-adic modular forms of
integral weight considered in \cite{BGK}. As explained in [\emph{loc.cit.}, p.~2394],
these can be directly related to the $p$-adic modular forms introduced by Serre \cite{Serre350}.
\end{remark}

%\end{comment}

Let $W_1 = X_{(p^{-p/p+1})}$ and $W_2 = X_{(p^{-1/p+1})}$, both open neighborhoods of $X^{\mathrm{ord}}$ with $W_2 \subseteq W_1$. Let
\begin{align*}
U: H^0(W_2, \underline{\omega}^k) &\longrightarrow H^0(W_1, \underline{\omega}^k) \subseteq H^0(W_2, \underline{\omega}^k) \\
V: H^0(W_1, \underline{\omega}^k) &\longrightarrow H^0(W_2, \underline{\omega}^k)
\end{align*}
be the operators defined in the introduction. Let $f\in S_k(\Gamma_1(N),K)$ be a newform defined over a number field $K$, 
%and of level $N$ coprime to $p$.
%Let $K_p$ be the $p$-adic completion of $K$, 
and consider $f$ as an element of $H^0(W_1,\underline{\omega}^k)$ by restriction. Then
$$
T_p(f) = U(f) + \chi(p)p^{k-1}V(f) \in H^0(W_2, \underline{\omega}^k),
$$
where $T_p$ is the $p$-th Hecke operator. In particular, if $f$ is an eigenform of level $\Gamma_0(N)$ and character $\chi$  with $T_p$-eigenvalue equal to $a_p$ then
$$
a_p\, f = U(f) + \chi(p)p^{k-1}V(f) \in H^0(W_2, \underline{\omega}^k).
$$

%We say that an overconvergent modular form in $H^0(W_2, \underline{\omega}^k)$ has {\em slope \rho} if it is an eigenvector of $U$ with an eigenvalue of %$p$-adic valuation $\rho$.

\begin{proposition}\label{prop:p-stabilize}
Let $f = \sum_{n=1}^\infty a_nq^n\in S_k(\Gamma_1(N),K)$ be a newform and let
\[
T^2 -a_pT + \chi(p)p^{k-1} = (T-\beta)(T-\beta')
\]
be the $p$-th  Hecke polynomial of $f$. %Assume that $|\alpha| < |\beta|$.
Then the overconvergent modular forms
$$
f_{\beta} := f - \beta' V(f), \quad
f_{\beta'} := f - \beta V(f)
$$
in $H^0(W_2,\underline{\omega}^k)$ are $U$-eigenvectors with eigenvalues $\beta$ and $\beta'$, respectively.
\end{proposition}

\begin{proof}
This follows from a straightforward calculation. Indeed,
viewing $f$ and $V(f)$ as sections in $H^0(W_2,\underline{\omega}^k)$, we see that
\begin{align*}
Uf_\beta=Uf-\beta' UV(f)&=Uf-\beta' f\\
&=T_pf-\chi(p)p^{k-1}V(f)-\beta' f\\
&=(a_p-\beta')f-\chi(p)p^{k-1}V(f)\\
&=\beta f_{\beta},
\end{align*}
using the relations $a_p=\beta+\beta'$ and $\chi(p)p^{k-1}=\beta\beta'$ for the last equality. The proof for $f_{\beta'}$ is obviously the same.
\end{proof}

Let now $W = X_{(\rr)}$, with $0<\rr<1$, be an open neighborhood of $X^{\rm{ord}}$,
and for any $r\geqslant 0$ consider the space
$$
\mathbb{H}^1(W, \nabla_r):= \mathbb{H}^1( \cH_r|_{W} \stackrel{\nabla_r}\longrightarrow \cH_r|_{W}\otimes \Omega^1_{W}(\log C)).
$$

\begin{theorem}[See \cite{coleman2}, \S 5]\label{theorem:p-adicCohomology}
\noindent

\begin{itemize}
\item[(i)]
There is a canonical isomorphism
$$
\mathbb{H}^1(W, \nabla_r) \simeq \frac{H^0(W, \underline{\omega}^{r+2})}{\theta^{r+1} H^0(W, \underline{\omega}^{-r}) }.
$$

\item[(ii)]  For any two open neighborhoods $W, W'$ of $X^{\rm{ord}}$, there is a canonical isomorphism
$$
\mathbb{H}^1(W, \nabla_r) = \mathbb{H}^1(W', \nabla_r).
$$
\end{itemize}
\end{theorem}

By restriction, there is an injection
$$
\mathbb{H}_{\rm{par}}^1(X,\nabla_r)\hookrightarrow \mathbb{H}^1(W\smallsetminus C, \nabla_r)
$$
for any choice of open neighborhood $W$ of $X^{\rm ord}$. The image of this map can be characterized by $p$-adic residues
(\cite[Prop.~3.9]{bdp}). In particular, if $f\in S_k(\Gamma_1(N),K)$ is a newform of weight $k\geqslant 2$,
the cohomology classes
$$
\{[f_{\beta}], [f_{\beta'}]\} \subseteq \mathbb{H}^1(W_2\smallsetminus C,\nabla_{k-2})
$$
naturally lie in $\mathbb{H}_{\rm{par}}^1(X_{K},\nabla_{k-2})$, and more precisely 
they lie in the $f$-isotypical component $\mathbb{H}_{\rm{par}}^1(X_{K},\nabla_{k-2})_f$,
%(\cite[Prop.~3.9]{bdp}),
which is a two-dimensional $K$-vector space.

\begin{theorem}
\label{theorem:canonicalBasis}
Let $f = \sum_{n=1}^\infty a_nq^n\in S_k(\Gamma_1(N),K)$ be a newform of weight $k\geqslant 2$,
and let $\beta$ and $\beta'$ be the roots of $T^2 -a_pT + \chi(p)p^{k-1}$,
ordered so that $v_p(\beta)\leqslant v_p(\beta')$.
Assume that the following two conditions hold:
\begin{itemize}
\item[(i)] $\beta\neq\beta'$.
\item[(ii)] $v_p(\beta')\neq k-1$.
\end{itemize}
Then $\{[f], [V(f)]\}$ is a basis for $\mathbb{H}_{\rm{par}}^1(X_K,\nabla_{k-2})_f$.
\end{theorem}

\begin{proof}
Since $\mathbb{H}_{\rm{par}}^1(X_K,\nabla_{k-2})_f$ is two-dimensional, it suffices to show that 
$[f]$ and $[V(f)]$ are linearly independent.
By Proposition~\ref{prop:p-stabilize} and \cite[Lem.~6.3]{coleman2}, condition (ii) guarantees that $[f_{\beta'}]\neq 0$,
and therefore it is an eigenvector of $Ver$ (see [\emph{loc.cit.}, Thm.~5.4]) acting on parabolic cohomology with eigenvalue $\beta'$.
In the same manner, the class $[f_{\beta}]$ is non-trivial, and it is an eigenvector of $Ver$ with eigenvalue $\beta$.
Thus by condition (i), the classes $[f_\beta]$ and $[f_{\beta'}]$ are linearly independent, and so must be
$[f]$ and $[V(f)]$, since $f_\beta$ and $f_{\beta'}$ are linear combinations of $f$ and $V(f)$.
\end{proof}

\begin{remark}\label{rem:variational-Hodge}
As clear from the proof, condition (ii) in Theorem~\ref{theorem:canonicalBasis}
could be weakened to the following:
\begin{itemize}
\item[(ii')] either $v_p(\beta')\neq k-1$ or $[f_{\beta'}]\neq 0$.
\end{itemize}
By \cite[Prop.~7.1]{coleman2}, condition (ii') fails if $f$ has CM by an imaginary quadratic field in which 
$p$ splits; conjecturally (see e.g. \cite{emerton-note}), these are the \emph{only} cases in which condition (ii') fails, 
but this is not known in general.
\end{remark}

\section{Recovering the shadow}
\label{section:recovering the shadow}

Let $f=\sum_{n=1}^{\infty} a(n)q^n$ 
be a newform satisfying the hypotheses of Theorem \ref{theorem:canonicalBasis},
and let $F$ be a harmonic Maass form which is good for $f$, in the sense of Definition \ref{def:goodHWMF}. 
The shadow $f$ can be recovered from $F$ by
$$
\xi_{2-k}(F) = \frac{f}{\|f\|^2}
$$
where $\|f\|$ is the Petersson norm of $f$. By the results in Section~\ref{section:harmonicMaassForms},
the harmonic Maass form $F$ has a holomorphic part $F^+$ with the property that
\begin{equation}\label{eq:prop3}
D^{k-1}(F^{+}) = \phi -  s_1\, f
\end{equation}
for some $\phi\in S_k^{!}(\Gamma_1(N),K)$ and some $s_1\in\mathbb{C}$.

In \cite[Thm.~1.2]{GuerzhoyKentOno}, Bringmann, Guerzhoy and Kane prove that one of the two $p$-stabilizations of $f$
can be recovered $p$-adically from an iterated application of $U$ to a certain `correction' of $D^{k-1}(F^+)$. 
In this section, we deduce their result from the $p$-adic techniques developed above. 
We begin by giving a new proof of \cite[Thm.~1.1]{GuerzhoyKentOno}.

%\footnote{This is perhaps not so surprising, as the original proof of \cite[Thm.~1.1]{GuerzhoyKentOno}
%follows from a modification of the proof of \cite[Thm.~1.3]{BruinierOnoRhoades}.}

\begin{theorem}\label{thm:GKO-1}
Let $\alpha\in\bC$ be such that $\alpha-c^+(1)\in K$.
Then the coefficients of
\[
\mathcal{F}_\alpha:=F^+-\alpha E_f:=\sum_{n\gg -\infty}c^+(n)q^n-\alpha\sum_{n=1}^\infty a(n)n^{1-k}q^n
\]
are all in $K$.
\end{theorem}

\begin{proof}
Write $\phi = \sum_{n\gg -\infty} d(n)q^n$, with $d(n)\in K$. By \eqref{eq:prop3}, we have the formula
\begin{equation}
c^+(n)=\left(\frac{d(n) - s_1 a(n)}{n^{k-1}}\right)\label{48}
\end{equation}
where $a(n):=0$ for $n\leqslant 0$.
The result is thus clear for $n\leqslant 0$.
Now let $n\geqslant 1$, and write $\alpha=c^+(1)+\gamma$ with $\gamma\in K$, or equivalently,
$\alpha=d(1) - s_1 +\gamma$.
Using $(\ref{48})$, an immediate calculation then reveals that the coefficient of
$q^n$ in $\mathcal{F}_\alpha$ is given by
\[
\frac{d(n)-d(1)-\gamma}{n^{k-1}}
\]
and the result follows.
\end{proof}

%We note that the function
%\[
%E_f = \sum_{n=1}^{\infty} a(n) n^{1-k}q^n
%\]
%appearing in the statement of Theorem~\ref{thm:GKO-1} is the so-called the Eichler integral of $f$,
%and it satisfies $D^{k-1}(E_f)=f$ by construction.

Since one can always take $\alpha=c^+(1)$ in Theorem~\ref{thm:GKO-1},
the coefficients of $\mathcal{F}_{c^+(1)}$ are all in $K$,
and so they may be viewed in $\bC_p$ via our fixed embedding $\overline{\mathbb{Q}}\hookrightarrow\bC_p$.

The following result is a special case of \cite[Thm.~1.2]{GuerzhoyKentOno}, but the ideas
in the proof will allow us to recover their result in its full strength (see Theorem~\ref{thm:GKO-2} below).

\begin{theorem}\label{shadow}
Assume that $v_p(\beta)<v_p(\beta')$ and that $v_p(\beta')\neq k-1$. Then
\[
\lim_{w\to +\infty}\frac{U^wD^{k-1}(\mathcal{F}_{c^+(1)})}{c_{c^+(1)}(p^w)} = f_{\beta},
\]
where we write $D^{k-1}(\mathcal{F}_{c^+(1)})= \sum_{n\gg-\infty} c_{c^+(1)}(n)q^n$.
\end{theorem}

\begin{proof}
First note that by equation $(\ref{eq:prop3})$ and $(\ref{48})$, we have
\[
D^{k-1}(\mathcal{F}_{c^+(1)}) = \phi - d(1) f,
\]
which is a weakly holomorphic cusp form of weight $k$ with coefficients in $K$,
defining a class in $\mathbb{H}_{\rm{par}}^1(X_K,\nabla_{k-2})_f$.
Now our assumptions clearly imply conditions (i) and (ii) of Theorem \ref{theorem:canonicalBasis},
and so (as shown in the proof of that result) the space $\mathbb{H}_{\rm{par}}^1(X_K,\nabla_{k-2})_f$ 
has a basis $\{ [f_{\beta}], [f_{\beta'}] \}$ of eigenvectors for $U$.
In particular, we can write
$$
[D^{k-1}(\mathcal{F}_{c^+(1)})] = t_1 [f_{\beta}] + t_2[f_{\beta'}] %\in \mathbb{H}_{\rm{par}}^1(X_K,\nabla_{k-2})_f.
$$
for some constants $t_1, t_2\in K$. The differential
$D^{k-1}(\mathcal{F}_{c^+(1)}) - t_1 f_{\beta} - t_2f_{\beta'}$ defines a class in $\HH^1(W_2\smallsetminus C, \nabla_{k-2}) = H^0(W_2\smallsetminus C,\underline{\omega}^k)/\theta^{k-1}H^0(W_2\smallsetminus C,\underline{\omega}^{2-k})$. This class is exact, by construction, and thus we may write
\begin{equation}
\label{equation:pAdicPhi}
D^{k-1}(\mathcal{F}_{c^+(1)}) = t_1f_{\beta} + t_2f_{\beta'} + \theta^{k-1}h
\end{equation}
for some  $h \in H^0(W_2\smallsetminus C,\underline{\omega}^{2-k})$. Applying $U$ to both sides of the equation gives
$$
UD^{k-1}(\mathcal{F}_{c^+(1)}) = t_1\beta f_{\beta} + t_2 \beta' f_{\beta'} + U(\theta^{k-1}h)
$$
and more generally, for any power $w\geqslant 1$, we obtain
$$
U^wD^{k-1}(\mathcal{F}_{c^+(1)}) = t_1 \beta^w f_{\beta} + t_2\beta'^w f_{\beta'} + U^w(\theta^{k-1}h).
$$
Dividing by $\beta^{w}$ we get
$$
\beta^{-w} U^wD^{k-1}(\mathcal{F}_{c^+(1)}) = t_1 f_{\beta} + t_2 \left(\frac{\beta'}{\beta}\right)^w f_{\beta'} + \beta^{-w}U^w(\theta^{k-1}h)
$$
and taking the limit as $w\to+\infty$ gives
$$
\lim_{w\to+\infty} \beta^{-w} U^wD^{k-1}(\mathcal{F}_{c^+(1)}) = t_1 f_{\beta}.
$$
This is because $v_p(\beta'/\beta)>0$ by the hypotheses and the differential $U^w(\theta^{k-1}h)$ 
has bounded denominators but its coefficients have arbitrarily high valuation as $w\to+\infty$.

To determine the value of the constant $t_1$, consider the coefficient of $q^p$ in  \eqref{equation:pAdicPhi}, which is
given by
\begin{align*}
c_{c^+(1)}(p) &= t_1 (a_p - \beta') + t_2(a_p - \beta) + O(p^{k-1}) \\
&= t_1 \beta + t_2 \beta' + O(p^{k-1}).
\end{align*}
By applying the multiplicative properties of the Fourier coefficients of newforms we get
$$
c_{c^+(1)}(p^w) = t_1\beta^w + t_2 \beta'^w + O(p^{w(k-1)})
$$
and taking the limit we obtain
$$
\lim_{w\to+\infty}\beta^{-w}c_{c^+(1)}(p^w) = t_1
$$
which gives the result.
\end{proof}

%\begin{remark}\label{non-critical}
%The assumption $v_p(\beta')\neq k-1$ in our Theorem~\ref{shadow} 
%seems to be missing in the statement of the corresponding [\cite[Thm.~1.2]{GuerzhoyKentOno}].
%\end{remark}

%\section{Overconvergent $p$-harmonic Maass forms}

%Let $f \in S_k(\Gamma_0(N),K)$ be a newform. Let $\phi\in S^{!}_k(\Gamma_0(N))$ be a weakly holomorphic cusp form and suppose that its class
%$$
%[\phi] \in \mathbb{H}_{\rm{par}}^1(X,\nabla_k)_f
%$$
%in parabolic cohomology lies in the $f$-isotypical component. Suppose additionally that $f$ satisfies the hypotheses of Theorem \ref{theorem:canonicalBasis}. %We may then write
%$$
%[\phi] = [af + b V(f)] \in \mathbb{H}_{\rm{par}}^1(X,\nabla_k)_f,
%$$
%for some $a,b\in K$. Let now $$W_2 - C = Y_{(p^{-1/p+1})}.$$ The class
%$$
%[\phi-af - b V(f)] \in \mathbb{H}^1(Y_{(p^{-1/p+1})},\nabla_k) = H^0(Y_{(p^{-1/p+1})}, \underline{\omega}^{k})/H^0(Y_{(p^{-1/p+1})}, %\theta^{k-1}\underline{\omega}^{2-k}),
%$$
%represented by the rigid analytic differential $\phi - af - bV(f)$, is exact, and we may thus find an element $F\in %H^0(Y_{(p^{-1/p+1})},\underline{\omega}^{2-k})$ such that
%$$
%\theta^{k-1} F = \phi - a\,f - b\,V(f)  \in H^0(Y_{(p^{-1/p+1})}, \underline{\omega}^{2-k}).
%$$
%\begin{definition}
%The weakly overconvergent modular form $F$ is a {\em $p$-harmonic Maass form} corresponding to the cusp form $f$.
%\end{definition}

%DISCUSS $a,b$.

Now we modify slightly the argument in Theorem~6.1 to recover \cite[Thm.~1.2]{GuerzhoyKentOno} in its full strength.
This refinement will be key for the results relating mock modular forms to $p$-adic modular forms in the next section.

%Throughout the following, we assume that $p\nmid N$, and that the roots $\beta$ and $\beta'$ of the $p$-th Hecke polynomial of
%$f\in S_k(\Gamma_0(N))$ are such that $v_p(\beta)<v_p(\beta')$, or equivalently, that $v_p(\beta)\neq (k-1)/2$.

For any $\alpha$ with $\alpha-c^+(1)\in K$, define
\[
\mathcal{F}_\alpha:=F^+-\alpha E_f
\]
and let $c_\alpha(n)$ denote the $n$-th coefficient in the expansion
\[
D^{k-1}(\mathcal{F}_\alpha)=\sum_{n\gg-\infty}c_\alpha(n)q^n.
\]

\begin{theorem}\label{thm:GKO-2}
Assume that $v_p(\beta)\neq v_p(\beta')$ and that $v_p(\beta')\neq k-1$.
Then for all but at most one choice of $\alpha$ with $\alpha-c^+(1)\in K$, we have
\[
\lim_{w\to+\infty}\frac{U^wD^{k-1}(\mathcal{F}_\alpha)}{c_\alpha(p^w)}=f_{\beta}.
\]
\end{theorem}

\begin{proof}
As in the proof of Theorem~6.1, we can write
\begin{equation}\label{eq:lc}
[D^{k-1}(\mathcal{F}_{c^+(1)})]=t_1[f_\beta]+t_2[f_{\beta'}]\in\mathbb{H}^1_{\rm par}(X_K,\nabla_{k-2})_f
\end{equation}
with
\[
t_1=\lim_{w\to+\infty}\frac{c_{c^+(1)}(p^w)}{\beta^w}.
\]
%where $c_{c^+(1)}(n)$ denotes the $n$-th coefficient of $D^{k-1}(\mathcal{F}_{c^+(1)})$. %=D^{k-1}(\mathcal{F}^+)-c^+(1)f$.
Let $\gamma\in K$ be such that $\alpha=c^+(1)+\gamma$, so that
$\mathcal{F}_\alpha=\mathcal{F}_{c^+(1)}-\gamma E_f$ by definition.
Noting that
\[
f=\frac{\beta f_{\beta}-\beta'f_{\beta'}}{\beta-\beta'},
\]
and substituting into the expression corresponding to $(\ref{eq:lc})$ for $\mathcal{F}_\alpha$
in place of $\mathcal{F}_{c^+(1)}$, we obtain
\[
[D^{k-1}(\mathcal{F}_\alpha)]=\left(t_1-\gamma\frac{\beta}{\beta-\beta'}\right)[f_\beta]
+\left(t_2+\gamma\frac{\beta'}{\beta-\beta'}\right)[f_{\beta'}],
\]
and hence we have the equality
\begin{equation}\label{eq:diff}
D^{k-1}(\mathcal{F}_\alpha)=\left(t_1-\gamma\frac{\beta}{\beta-\beta'}\right)f_\beta
+\left(t_2+\gamma\frac{\beta'}{\beta-\beta'}\right)f_{\beta'}+\theta^{k-1}h
\end{equation}
as sections in $H^0(W_2\smallsetminus C,\underline{\omega}^k)$, for some $h\in H^0(W_2\smallsetminus C,\underline{\omega}^{2-k})$.
Applying $U^w$ to both sides of this equation and letting $w\to+\infty$ as in the proof of Theorem~6.1,
we deduce that
\begin{equation}\label{eq:1}
\lim_{w\to+\infty}\frac{U^wD^{k-1}(\mathcal{F}_\alpha)}{\beta^w}=\left(t_1-\gamma\frac{\beta}{\beta-\beta'}\right)f_{\beta}.
\end{equation}
On the other hand, arguing again as in Theorem~6.1 we find that the $p^w$-th coefficient of $D^{k-1}(\mathcal{F}_\alpha)$
is given by
\[
c_\alpha(p^w)=\left(t_1-\gamma\frac{\beta}{\beta-\beta'}\right)\beta^w
+\left(t_2+\gamma\frac{\beta'}{\beta-\beta'}\right)\beta'^w+O(p^{w(k-1)}),
\]
and hence
\begin{equation}\label{eq:2}
\left(t_1-\gamma\frac{\beta}{\beta-\beta'}\right)=\lim_{w\to+\infty}\frac{c_\alpha(p^w)}{\beta^w}.
\end{equation}
Therefore, \emph{except} in the case where
\begin{equation}\label{eq:exception-gamma}
\gamma=\frac{t_1(\beta-\beta')}{\beta}=(\beta-\beta')\lim_{w\to+\infty}\frac{c_{c^+(1)}(p^w)}{\beta^{w+1}},
\end{equation}
combining $(\ref{eq:1})$ and $(\ref{eq:2})$ we recover $f_{\beta}$
from $\mathcal{F}_\alpha$ as in the statement of the theorem.
\end{proof}

\section{Mock modular forms as $p$-adic modular forms}
\label{section:mock modular forms as }

We now let $\alpha$ range over the larger set of values
\[
c^+(1)+\bC_p:=\{c^+(1)+\gamma\;\colon\;\gamma\in\bC_p\},
\]
and interpret the exceptional value of $\alpha$ in Theorem~\ref{thm:GKO-2} as the only one for which
the `corrected' mock modular form
\[
\mathcal{F}_\alpha=F^+-\alpha E_f
\]
%(which has coefficients in $K$ for any $\alpha\in c^+(1)+K$ by Theorem~\ref{thm:GKO-1}) 
gives rise to a $p$-adic modular form upon $p$-stabilization. 
%However, since we assume that $p\nmid N$, for that we will need to first take
%a certain $p$-stabilization of $\mathcal{F}_\alpha$. 
Recall that we let $\beta$ and $\beta'$ be the roots of the $p$-th Hecke polynomial of $f$,
ordered so that $v_p(\beta)\leqslant v_p(\beta')$.

\begin{definition}
For any $\alpha\in c^+(1)+\bC_p$, define
\[
\mathcal{F}_\alpha^*:=\mathcal{F}_\alpha-p^{1-k}\beta'\mathcal{F}_\alpha\vert V
\]
and write
\[
D^{k-1}(\mathcal{F}_\alpha^*)=\sum_{n\gg-\infty}c_\alpha^*(n)q^n.
\]
\end{definition}

Our first result shows that, similarly as in Theorem~\ref{thm:GKO-2} for $\mathcal{F}_\alpha$, 
the $p$-stabilization $f_{\beta}$ of the shadow of $F^+$
can be recovered $p$-adically from $\mathcal{F}_\alpha^*$.

\begin{theorem}\label{thm:GKO-2*}
Assume that $v_p(\beta)\neq v_p(\beta')$ and that $v_p(\beta')\neq k-1$.
Then for all but at most one choice of $\alpha\in c^+(1)+\bC_p$, we have
\[
\lim_{w\to+\infty}\frac{U^wD^{k-1}(\mathcal{F}_\alpha^*)}{c_\alpha^*(p^w)}=f_{\beta}.
\]
\end{theorem}

\begin{proof}
The proof is quite similar to the proof of Theorem~\ref{thm:GKO-2}.
Writing $\alpha=c^+(1)+\gamma$ with $\gamma\in\bC_p$, an immediate calculation reveals that
\begin{equation}\label{eq:F}
D^{k-1}(\mathcal{F}_\alpha^*)=D^{k-1}(\mathcal{F}_{c^+(1)})\vert(1-\beta'V)-\gamma f_{\beta}.
\end{equation}
As in the proof of Theorem~\ref{shadow}, we write
\[
[D^{k-1}(\mathcal{F}_{c^+(1)})]=t_1[f_\beta]+t_2[f_{\beta'}]\in\mathbb{H}^1_{\rm par}(X_K,\nabla_{k-2})_f,
\]
with $t_1=\lim_{w\to+\infty}\beta^{-w}c_{c^+(1)}$.
Applying the operator  $1-\beta'V$ to this last equality, and noting that $V=U^{-1}$ on cohomology, we obtain
\[
[D^{k-1}(\mathcal{F}_{c^+(1)})\vert(1-\beta'V)]=t_1\frac{(\beta-\beta')}{\beta}[f_{\beta}],
\]
and hence by $(\ref{eq:F})$:
\begin{equation}\label{eq:0}
[D^{k-1}(\mathcal{F}_\alpha^*)]=\left(\frac{t_1(\beta-\beta')}{\beta}-\gamma\right)[f_{\beta}].
\end{equation}
Arguing again as in the proof of Theorem~6.1, we obtain the equalities
\begin{equation}\label{eq:1*}
\lim_{w\to+\infty}\frac{U^w(D^{k-1}(\mathcal{F}_\alpha^*))}{\beta^w}
=\left(\frac{t_1(\beta-\beta')}{\beta}-\gamma\right) f_{\beta}
\end{equation}
and
\begin{equation}\label{eq:2*}
\frac{t_1(\beta-\beta')}{\beta}-\gamma=\lim_{w\to+\infty}\frac{c_{\alpha}^*(p^w)}{\beta^w}.
\end{equation}
Therefore, \emph{except} in the case where
\begin{equation}\label{eq:exception}
\gamma=\frac{t_1(\beta-\beta')}{\beta}=(\beta-\beta')\lim_{w\to+\infty}\frac{c_{c^+(1)}(p^w)}{\beta^{w+1}},
\end{equation}
the combination of $(\ref{eq:1})$ and $(\ref{eq:2})$ recovers $f_{\beta}$
from $\mathcal{F}_\alpha^*$ as in the statement of the theorem.
\end{proof}

Considering the exceptional value of $\alpha$ arising in the proof of Theorem~\ref{thm:GKO-2*}, 
we recover the result of \cite[Thm.~1.1]{BGK}. %We continue to assume that $v_p(\beta)<v_p(\beta')$.

\begin{theorem}\label{thm:BGK-1}
Assume that $v_p(\beta)\neq v_p(\beta')$ and that $v_p(\beta')\neq k-1$.
Then among all values of $\alpha\in c^+(1)+\bC_p$, the value
\[
\alpha=c^+(1)+(\beta-\beta')\lim_{w\to+\infty}\frac{c_{c^+(1)}(p^w)}{\beta^{w+1}}
\]
is the unique one such that $\mathcal{F}_\alpha^*$ is a $p$-adic
modular form
%\footnote{we should explain at some point at our $p$-adic modular forms are allowed
%to have poles at the cups, as in the first Remark in the introduction to \cite{BGK}.}
of weight $2-k$.
\end{theorem}

\begin{proof}
Write $\alpha=c^+(1)+\gamma$ with $\gamma\in\bC_p$.
Since $[f_{\beta}]\neq 0$ (see the proof of Theorem~\ref{theorem:canonicalBasis}),
we deduce from $(\ref{eq:0})$ and $(\ref{eq:exception})$ that the class of $D^{k-1}(\mathcal{F}_\alpha^*)$
in $\mathbb{H}^1_{\rm par}(X_K,\nabla_{k-2})$ vanishes only for the value of $\alpha$ given in the statement.
Now, since the natural restriction map
\[
\mathbb{H}^1_{\rm par}(X_K,\nabla_k)\longrightarrow\mathbb{H}^1(W_2\smallsetminus C,\nabla_{k-2})
=\frac{H^0(W_2\smallsetminus C,\underline{\omega}^k)}{\theta^{k-1}H^0(W_2\smallsetminus C,\underline{\omega}^{2-k})}
\]
is injective, the above value of $\alpha$ is also the unique one such that the class of
$D^{k-1}(\mathcal{F}_\alpha^*)$ becomes trivial in $\mathbb{H}^1(W_2\smallsetminus C,\nabla_{k-2})$,
and hence such that $\mathcal{F}_\alpha^*\in H^0(W_2\smallsetminus C,\underline{\omega}^{2-k})$.
%as was to be shown.
\end{proof}

Next we consider a second modification of $\mathcal{F}_\alpha$.

\begin{definition}\label{def:F-alpha-delta}
For any $\delta\in\bC_p$, define
\[
\mathcal{F}_{\alpha,\delta}:=\mathcal{F}_\alpha-\delta(E_f-\beta E_{f\vert V}).
\]
\end{definition}

Our next result explores the values of $\alpha$ and $\delta$ for which
$\mathcal{F}_{\alpha,\delta}$ is a $p$-adic modular form, recovering the content of \cite[Thm~1.2(2)]{BGK}.

\begin{theorem}\label{thm:BGK-2}
Assume that $v_p(\beta)\neq v_p(\beta')$ and that $v_p(\beta')\neq k-1$. Then there exists a unique pair of values $(\alpha,\delta)$
for which $\mathcal{F}_{\alpha,\delta}$ is a $p$-adic modular. In fact, $\alpha$ is as in Theorem~\ref{thm:BGK-1},
and
\[
\delta=\lim_{w\to+\infty}\frac{a_{\mathcal{F}_\alpha}(p^w)p^{w(k-1)}}{\beta'^w}.
\]
Here, we write $\mathcal{F}_\alpha=\sum_{n\gg -\infty}a_{\mathcal{F}_\alpha}(n)q^n$.
\end{theorem}

\begin{proof}
With the same notations as in the proof of Theorem~\ref{thm:GKO-2}, we can write
the equality
\begin{equation}\label{eq:lc-2}
[D^{k-1}(\mathcal{F}_{\alpha,\delta})]=\left(t_1-\gamma\frac{\beta'}{\beta-\beta'}\right)[f_\beta]
+\left(t_2+\gamma\frac{\beta}{\beta-\beta'}-\delta\right)[f_{\beta'}]
\end{equation}
in $\mathbb{H}^1_{\rm par}(X_K,\nabla_{k-2})_f$. Since we may check the triviality of these classes
upon restriction to $W_2\smallsetminus C$, it follows that $\mathcal{F}_{\alpha,\delta}$ is a $p$-adic modular form of weight $2-k$
if and only if the class $[D^{k-1}(\mathcal{F}_{\alpha,\delta})]$ vanishes. As in the proof of Theorem~\ref{theorem:canonicalBasis}, %under our running assumptions
the classes $[f_\beta], [f_{\beta'}]$ form a basis for $\mathbb{H}^1_{\rm par}(X_K,\nabla_{k-2})_f$,
and hence $\mathcal{F}_{\alpha,\delta}$ is a $p$-adic modular form if and only if the coefficients in the
right-hand side of $(\ref{eq:lc-2})$ both vanish. In particular (second coefficient), this shows that the value of $\gamma$
is given by $(\ref{eq:exception-gamma})$, and therefore the necessary value of
$\alpha=c^+(1)+\gamma$ is the same as in Theorem~\ref{thm:BGK-1}. To determine the value of $\delta$, we first rewrite Equation $(\ref{eq:diff})$ for the above value of $\alpha$
(so that the first summand in the right-hand side of that equation vanishes):
\[
D^{k-1}(\mathcal{F}_\alpha)=\left(t_2+\gamma\frac{\beta'}{\beta-\beta'}\right)f_{\beta'}+\theta^{k-1}h.
\]
Equating the $p^w$-th coefficients in this equality, we obtain
\[
c_\alpha(p^w)=\left(t_2+\gamma\frac{\beta'}{\beta-\beta'}\right)\beta'^w+O(p^{w(k-1)})
\]
and hence dividing by $\beta'^w$ and letting $w\to+\infty$ we deduce
\begin{equation}\label{eq:first}
\lim_{w\to+\infty}\frac{c_\alpha(p^w)}{\beta'^w}=\left(t_2+\gamma\frac{\beta'}{\beta-\beta'}\right).
\end{equation}
(Note that the assumption $v_p(\beta')<k-1$ is being used here.) Finally, substituting $(\ref{eq:first})$
into $(\ref{eq:lc-2})$, we see that the necessary value for $\delta$ is given by
\[
\delta=\lim_{w\to\infty}\frac{c_\alpha(p^w)}{\beta'^w}
=\lim_{w\to\infty}\frac{a_{\mathcal{F}_\alpha}(p^w)p^{w(k-1)}}{\beta'^w},
\]
as was to be shown.
\end{proof}

\section{The CM case}
\label{section:CM}

In this section we treat the case in which $f$ has CM. This case is of special interest,
since then one can choose a good harmonic Mass form $F$ for $f$ as in
Section~\ref{section:harmonicMaassForms} with $F^+$ having algebraic coefficients.
Conjecturally, this characterize the CM property of $f$ 
(see \cite[p.6170]{GuerzhoyKentOno}). Thus assume that $f=\sum_{n=1}^\infty a_nq^n\in S_k(\Gamma_1(N),K)$ has CM by
an imaginary quadratic field $M$ of discriminant prime to $p$,
and let $F=F^++F^-$ be a good harmonic Maass form attached to $f$.
We also assume (upon enlarging $K$, if necessary) that $K$ contains a primitive $m$-th root of unity,
where $m=N\cdot{\rm disc}(M)$. Then by \cite[Thm.~1.3]{BruinierOnoRhoades},
$F^+$ has coefficients in $K$, and so $D^{k-1}(F^+)$ defines a class in
$\mathbb{H}^1_{\rm par}(X_K,\nabla_{k-2})_f$.
%which by
%the proof of Theorem~\ref{theorem:canonicalBasis}
%admits (under some assumptions) a basis given by the classes $[f_\beta]$ and $[f_{\beta'}]$.

We first treat the case in which $p$ is inert in $M$. In this case $a_p=\beta+\beta'=0$,
and so by the proof of Theorem~\ref{theorem:canonicalBasis},
the space $\mathbb{H}^1_{\rm par}(X_K,\nabla_{k-2})_f$
admits a basis given by the classes $[f_\beta]$ and $[f_{\beta'}]$. 

\begin{lemma}\label{lem:-t1+t2}
Assume that $p$ is inert in $M$, and write
$[D^{k-1}(F^+)]=t_1[f_\beta]+t_2[f_{\beta'}]$.
Then
\[
\lim_{w\to+\infty}\frac{a_{D^{k-1}(F^+)}(p^{2w+1})}{\beta^{2w+1}}=t_1-t_2.
\]
\end{lemma}

\begin{proof}
The proof will be obtained by arguments similar to the proof of Theorem~\ref{shadow},
but some adjustments are necessary due to the fact the condition $v_p(\beta)\neq v_p(\beta')$
is not satisfied in this case. Instead, we shall exploit the extra symmetry $\beta'=-\beta$.

Upon restriction to $W_2\smallsetminus C$, we can write
\begin{equation}\label{eq:F+}
D^{k-1}(F^+)=t_1f_\beta+t_2f_{\beta'}+\theta^{k-1}h
\end{equation}
for some $h\in H^0(W_2\smallsetminus C,\underline{\omega}^{2-k})$.
Taking $p^{2w+1}$-st coefficients in this identity, we immediately obtain
\begin{align*}
a_{D^{k-1}(F^+)}(p^{2w+1})&=t_1\beta^{2w+1}+t_2\beta'^{2w+1}+O(p^{(2w+1)(k-1)})\\
&=(t_1-t_2)\beta^{2w+1}+O(p^{(2w+1)(k-1)}),
\end{align*}
and hence dividing by $\beta^{2w+1}$ and letting $w\to+\infty$ the result follows.
\end{proof}

\begin{definition}
For any $\alpha\in\bC_p$, define
\[
\widetilde{\mathcal{F}}_\alpha:=F^+-\alpha E_{f\vert V}.
\]
\end{definition}

Armed with Lemma~\ref{lem:-t1+t2}, in Corollary~\ref{cor:GKO} below we will determine the values of $\alpha$
for which $\widetilde{\mathcal{F}}_{\alpha}$ is a $p$-adic modular form, thus recovering
\cite[Thm.~1.3]{BGK}. This will be an immediate consequence of the following result.

\begin{theorem}\label{thm:BGK-3}
Assume that $p\nmid N$ is inert in $M$, and for any $\widetilde{\alpha}\in\mathbb{C}_p$ define
\[
G_{\widetilde{\alpha}}:=F^+-\widetilde{\alpha}(E_f-\beta E_{f\vert V}).
\]
Then there exists a unique value of $\widetilde{\alpha}$
such that $G_{\widetilde{\alpha}}$ is a $p$-adic modular form of weight $2-k$, and it is given by
\[
\widetilde{\alpha}=\lim_{w\to+\infty}\frac{a_{D^{k-1}(F^+)}(p^{2w+1})}{\beta^{2w+1}}.
\]
\end{theorem}

\begin{proof}
We will deduce this result by first determining the values of $\alpha$ and $\delta$
for which the form $\mathcal{F}_{\alpha,\delta}$ of Definition~\ref{def:F-alpha-delta}
is a $p$-adic modular form. Note that this case is not covered by Theorem~\ref{thm:BGK-2},
since the proof of that result relies crucially
on the assumption that $v_p(\beta)<v_p(\beta')$. Instead, we will exploit again the fact
that $\beta'=-\beta$. Since $[f_\beta]$ and $[f_{\beta'}]$ form a basis for $\mathbb{H}^1_{\rm par}(X_K,\nabla_{k-2})_f$, 
Equation (\ref{eq:lc-2}) for $[D^{k-1}(\mathcal{F}_{\alpha,\delta})]$ still applies,
and in this case it reduces  
(noting that we may set $\gamma=\alpha$ by the algebraicity of $c^+(1)$) to
\begin{equation}\label{eq:lc-3}
[D^{k-1}(\mathcal{F}_{\alpha,\delta})]=\left(t_1-\frac{\alpha}{2}\right)[f_\beta]
+\left(t_2-\frac{\alpha}{2}-\delta\right)[f_{\beta'}].
\end{equation}
By Theorem~\ref{theorem:canonicalBasis}, the classes $[f]$ and $[V(f)]$ form
a basis for the space $\mathbb{H}^1_{\rm par}(X_K,\nabla_{k-2})_f$, and rewriting $(\ref{eq:lc-3})$
in terms of them we arrive at
\begin{equation}\label{eq:lc-4}
[D^{k-1}(\mathcal{F}_{\alpha,\delta})]=(t_1+t_2-\alpha-\delta)[f]+\beta(t_1-t_2-\alpha-\delta)[V(f)].
\end{equation}
Now, $\mathcal{F}_{\alpha,\delta}$ is a $p$-adic modular form of weight $2-k$ if and only
if both coefficients in the right-hand side of Equation $(\ref{eq:lc-4})$ vanish;
in particular, we need to have
\begin{equation}\label{eq:alpha+delta}
\alpha+\delta=t_1-t_2=\lim_{w\to+\infty}\frac{a_{D^{k-1}(F^+)}(p^{2w+1})}{\beta^{2w+1}},
\end{equation}
where we used Lemma~\ref{lem:-t1+t2} for the second equality. The necessary vanishing
of $(\ref{eq:lc-4})$ also forces the vanishing of $t_2$ and hence from $(\ref{eq:lc-3})$
we deduce that $\delta=-\frac{\alpha}{2}$, or equivalently, $\alpha+\delta=\frac{\alpha}{2}$.
Finally, noting that
\begin{align*}
\mathcal{F}_{\alpha,\delta}%&=F^+-(\alpha+\delta)E_f+\delta\beta E_{f\vert V}\\
&=F^+-\frac{\alpha}{2}\left(E_f-\beta E_{f\vert V}\right)=G_{\frac{\alpha}{2}},
\end{align*}
we conclude from $(\ref{eq:alpha+delta})$ that $G_{\widetilde{\alpha}}$ is a $p$-adic modular form
if and only if $\widetilde{\alpha}$ is given by the $p$-adic limit in the statement.
\end{proof}

\begin{corollary}\label{cor:GKO}
Assume that $p\nmid N$ is inert in $M$. Then there exists a unique value of $\alpha$
such that $\widetilde{\mathcal{F}}_\alpha$ is a $p$-adic modular form of weight $2-k$, and it is given by
\[
\alpha=\lim_{w\to+\infty}\frac{a_{D^{k-1}(F^+)}(p^{2w+1})}{\beta^{2w}}.
\]
\end{corollary}

\begin{proof}
Comparing the definitions of $\widetilde{\mathcal{F}}_\alpha$ and $G_{\widetilde{\alpha}}$, we see that
\[
G_{\widetilde{\alpha}}=\widetilde{\mathcal{F}}_\alpha-\widetilde{\alpha}E_f,
\]
with $\alpha=\widetilde{\alpha}\beta$. Since $E_f$ is easily seen to be
%(and hence $\widetilde{\alpha}E_f$, for any $\widetilde{\alpha}\in\bC_p$)
a $p$-adic modular form under our hypotheses (see \cite[Prop.~4.2]{BGK}, which remains true in our case $p\nmid N$), 
the result follows from Theorem~\ref{thm:BGK-3}.
\end{proof}
We conclude this section by dealing with the case in which $f$ has CM
by an imaginary quadratic field $M$ in which $p$ splits, characterizing the values of $\alpha\in\bC_p$
for which $\mathcal{F}_{\alpha}^*$ is a $p$-adic modular form. As noted in Remark~\ref{rem:variational-Hodge}, 
the class $[f_{\beta'}]$
vanishes in this case, and so the proofs of Theorem~\ref{thm:GKO-2*} and Theorem~\ref{thm:BGK-1} break down.
However, based on the observation that (using the algebraicity of $c^+(1)$ to set $\alpha=\gamma$)
\begin{equation}
\label{eq:observ}
\mathcal{F}_\alpha^*=(F^+-\alpha E_f)\vert(1-p^{1-k}\beta'V)=\mathcal{F}_0^*-\alpha E_{f_{\beta}},
\end{equation}
%(since $c^+(1)$ is algebraic), 
we can easily prove the following result (cf. \cite[Thm.~1.2]{BGK}).

\begin{theorem}
Assume that $p\nmid N$ splits in $K$. Then among all values of $\alpha\in\bC_p$,
the value $\alpha=0$ is the unique one for which $\mathcal{F}_\alpha^*$ is a $p$-adic modular form
of weight $2-k$.
\end{theorem}

\begin{proof}
As we have already argued in preceding proofs, $\mathcal{F}_{\alpha}^*$ is a $p$-adic
modular form of weight $2-k$ if and only if the class $[D^{k-1}(\mathcal{F}_\alpha^*)]$
vanishes, %in $\mathbb{H}^1(W_2-C,\underline{\omega}^k)$,
and from $(\ref{eq:observ})$ we see that
\begin{equation}\label{eq:equiv}
[D^{k-1}(\mathcal{F}_\alpha^*)]=0\quad\Longleftrightarrow\quad\alpha[f_{\beta}]=[D^{k-1}(\mathcal{F}_0^*)].
\end{equation}
In particular, this shows that $\mathcal{F}_\alpha^*$ is a $p$-adic modular form of weight $2-k$ for $\alpha=0$, and so
$[D^{k-1}(\mathcal{F}_0^*)]=0$. On the other hand, since $[f_{\beta}]\neq 0$ (see the proof of Theorem~\ref{theorem:canonicalBasis}),
equivalence $(\ref{eq:equiv})$ shows that $[D^{k-1}(\mathcal{F}_\alpha^*)]\neq 0$ for $\alpha\neq 0$,
yielding the result.
\end{proof}

\renewcommand\refname{References}

\bibliographystyle{alpha}
\bibliography{p-harmonic}

\end{document}